\documentclass[twoside,12pt]{article}

\usepackage[hyphens]{url}
\usepackage{jmlr2e}
\usepackage{enumitem}
\usepackage{times,mathtools}
\usepackage{color}
\usepackage{options}
\usepackage{caption} 
\usepackage{subcaption} 
\usepackage{packages}
\usepackage{macros}
\usepackage{microtype,xparse,tcolorbox}
\usepackage{graphicx}
\usepackage{footnote}
\usepackage{ragged2e}
\usepackage{multirow}
\usepackage{float}
\usepackage{authblk}
\makeatletter
\newcommand*{\rom}[1]{\expandafter\@slowromancap\romannumeral #1@}

\allowdisplaybreaks

\begin{document}

\thispagestyle{plain}
\title{An Optimization Model for Scheduling Freight Trains on a Single Rail Track}

\newcommand{\HRule}{\rule{\linewidth}{0.5mm}} 
\center 

\HRule \\[0.4cm]
{\LARGE \bfseries An Optimization Model for Scheduling Freight Trains on a Single Rail Track}
\HRule \\[1 cm]

 {Marie Alaghband$^1$, Babak Farhang Moghaddam$^2$\\
  
  $^1$University of Central Florida, Orlando, Florida, USA \\
  \texttt{marie.alaghband@knights.ucf.edu}\\
  $^2$Institute of Management and Planning Studies, Tehran, Iran\\
  \texttt{farhang@imps.ac.ir}\\[1 cm]
  }

\center
\justify{
\begin{abstract}
In many countries, a rail network consists of single lines with sidings where interactions between trains occur (meet, pass). In this paper, we study two issues of these networks: first, the scheduling of freight trains in a single line corridor while ensuring safe interactions and second, the allocation of freight to freight trains regarding the release dates of freight, weight of freight, and weight capacity of the trains. Both of these issues must be addressed when examining real-world freight train scheduling problems. The objective functions of this study are the minimization of a train's travelling time, the allocation of freight to freight trains, and the reduction of tardiness of freight at destination. Both scheduling and allocation problems are presented using integer linear programming models. In addition, an integrated novel heuristic algorithm has been proposed to solve them. Computational results demo\-nstrated through a generated dataset show both model validation and efficiency of the heuristic algorithm. The heuristic algorithm has been designed to incorporate the practical operational railway rules with modest modification. Although its outputs slightly differ from the exact solutio\-ns, it can solve both models simultaneously in large scale problems.
\end{abstract}
}


\begin{flushleft}
\justify{
\section{Introduction}

Transportation is vital for various human activities \cite{blazewicz2019scheduling}. A major type of transport is rail transport. Several different problems must be resolved in order to achieve efficient rail transportation; these can be modelled and solved individually. Based on the survey of \cite{assad1980models}, the rail modeling problems can be categorized in the following groups: Institutional Background, Facilities Location, Yard and Terminal Models, Line Models, Rail Network Model, Blocking and Train Formation, Train Schedules and Timetables, and Car and Engine Distribution. In this study, train scheduling and timetabling, which is one of the most important categories, is addressed. Some other literature summaries in this area of research were published by \cite{haghani1987rail}, \cite{cordeau1998survey}, \cite{lusby2011railway}, and \cite{harrod2010operations}. Various scheduling problems have been modelled and described in detail \cite{meisel2019scheduling}, \cite{pellegrini2019scheduling}. The first work that sought to find an optimum solution for the train scheduling problem was started by \cite{szpigel1973optimal}. He developed a linear programming model with a branch and bound method to minimize the sum of travel times. \cite{higgins1996optimal} then considered the train scheduling problem on a single line track. They proposed a multi-objective mathematical programming model in a branch and bound procedure, where the objective was to minimize the deviation from scheduled arrival time and fuel consumption costs.

Although there has been vast research on passenger train scheduling, only a few researchers have considered freight train scheduling. As rail traffic in both passenger and freight trains on mixed-use rails is growing continuously \cite{talebian2018capacity}, one branch of the research conducted in the literature focused on railway systems that both passenger and freight trains use. \cite{godwin2007freight} addressed the problem of scheduling freight trains in a passenger rail network. They showed that freight train scheduling in a passenger rail network is NP- complete and developed a step-wise dispatching heuristic considering several objectives (i.e., percentage deviation of sum of travel times from lower bound, percentage standard mean tardiness, percentage tardy trains, percentage conditional mean tardiness, and percentage maximum tardiness). \cite{xu2015scheduling} reported the design of an improved switchable policy which is rooted in approach\-es by \cite{mu2013efficient}, with the analysis of possible delays caused by different path choices. Also, \cite{fu2017models} studied how changing the speed limits of different railway segments affected efficiency.  \cite{cacchiani2010scheduling} presented an integer linear programming formulation to address the same problem, with the objective of scheduling and assigning as many new freight trains as possible on railway networks. In their study, they only considered the constraint of freight train capacity. To study the same problem for mixed-use rail systems, \cite{zyngier2018process} developed a detailed scheduling model with a process systems approach. They also solved their model for a week long period and were able to yield fast solutions with significant improvements to solution times. \cite{murali2016modeling} proposed an expert tool to help train schedule planners determine proper routes and schedules for short time frames, and to manage the restricted track capacity available for train movements. \cite{rahimi2016optimization} developed a technique based on discrete‐event simulation and response surface methodology to model and then optimize the schedule of subway train travels. In a recent study, \cite{behiri2018urban} formulated the problem of the freight rail transport scheduling using MIP and proved the NP-hardness of the problem. They then proposed two heuristics based on dispatching rules and single-train-based decomposition and evaluated their models using a discrete-event simulation approach.

Another area of research studies problems associated with railway systems that only service freight trains. \cite{jaumard2012dynamic} presented research on these freight train scheduling problems. They identified more comprehensive constraints (i.e., travel and dwelling time, safety distance, segment conflict, and capacity) with the usage of mixed integer programming. One study dates back to 2014, when \cite{rahman2014freight} represented an integer programming formulation for the freight train scheduling problem in a single line corridor. Regarding safe interactions between trains as constraints, the objective was to minimize the arrival time of the last train at its destination. Finally, \cite{ke2015new} addressed the problem of freight train timetabling on a single-track railway system to minimize the train waiting times. They presented a new method that utilized both fixed-block signaling systems and fuzzy logic systems to address the problem of freight train timetabling.

In terms of objective function, \cite{kuo2010freight} proposed that the most common objective functions of freight train scheduling and timetabling are mini\-mizing deviation from the schedule, operating cost, train delay, and average travel time.

Using the Greedy Algorithm, \cite{sinha2016iterative} presented an iterative bi-level hierarchical approach for train scheduling based on the decentralized operat\-ional control concept in railway operations where they divided the entire railway network into a number of sub-networks connected at boundary stations, called interchange points.

In this paper, we address railway network systems exclusive to freight trains. Two main problems are addressed: scheduling freight trains in a single line corridor to minimize the total train’s travel time, and allocating freight to scheduled freight trains to maximize the allocation of freight and minimize the tardiness of freight to their destination.
Since both models have their own complexities, their combination will be complex too. Hence, a novel heuristic algorithm has been proposed that can simultaneously address train allocation and scheduling on a single line corridor. Next, the methodologies for both the scheduling and allocation problems have been provided in two subsections. In the next section, which is subdivided into two subsections, we provide the methodology for both the scheduling problem and the allocation problem.

\begin{table}[ht]
\caption{Subscripts and parameters}
\label{parameters_table}
\centering
\scriptsize
\resizebox{0.8\textwidth}{!}{
\begin{tabular}{p{1.5cm}p{7.5cm}}
\hline
Symbol & Definition \\
\hline
$S$ & =\{segment 1, segment 2, segment 3\}, set of all segments \\ [5pt]
$T$ & =\{departing train 1, departing train 2, departing train 3, returning train 1, returning train 2\}, set of trains \\ [5pt]
$t$ & =\{departing train 1, departing train 2, departing train 3\}, set of departing trains \\ [5pt]
$r$ & =\{returning train 1, returning train 2\}, set of returning trains \\ [5pt]
$pd$ & Set of segments containing departing trains\\ [5pt]
$pdr$ & Set of segments containing returning trains\\ [5pt]
$SP$ & Start station\\ [5pt]
$EP$ & End station \\ [5pt]
$\pi_t$ & Priority of departing trains \\ [5pt]
$\upsilon_s^t$ & Average speed of train t at segment $pd$ \\ [5pt]
$dw_{pd}^t$ & Dwelling time of train t at segment $pd$ \\ [5pt]
$mr_{pd}^t$ & Minimum time for train t to travel segment $pd$ \\  [5pt]
$ST_{pd}$ & Safety lag time at segment $pd$ \\ [5pt]
$Lo_{pd}^t$ & Loading time of train t at segment $pd$ \\ [5pt]
$Ul_{pd}^t$ & Unloading time of train t at segment $pd$ \\ [5pt]
$M$ & Sufficiently large constant \\ [5pt]
$n$ & Number of segments \\ [5pt]
\end{tabular}%
}
\end{table}

\section{Methodology}
In this section, we first formally describe the developed model to address the scheduling problem and then the model and formulation regarding freight allocation.

\begin{table}[b]
\caption{Decision variables of the scheduling model}
\label{Decision_variables}
\centering
\scriptsize
\resizebox{0.8\textwidth}{!}{%
\begin{tabular}{p{1.5cm}p{7.5cm}}
\hline
Symbol & Definition \\
\hline
$d_{pd}^t$ & Departure time of departing train $t$ from segment  $pd$ \\ [5pt]
$a_{pd}^t$ & Arrival time of departing train $t$  to segment  $pd$ 
\end{tabular}%
}
\end{table}

\subsection{The scheduling problem} \label{scheduling problem}
In this study, we considered freight trains travelling on a single line corridor. Thus, the start and end stations are located at the start and end of this corridor, respectively. Like real world corridors, our model is divided into segments operat\-ed by stations. At most, one pair of trains can cross and overtake one another at these stations. Similar to the study of \cite{xu2019efficient}, to simplify the problem, some necessary assumptions were made: the route of each train is fixed, and since one single line corridor is considered, only two categories of trains are identified—departing and returning trains. Also, the traveling time at each station was assumed to be zero.

Based on the traveling directions of departing and returning trains, each segment has two different names corresponding to departing and returning trains. Two trains can follow each other at a minimum distance depending on the leading train's speed. 

Rail transportation is used to for a variety of goods in the real world, and due to the differences between goods (such as value, expiration dates, etc), their importance will also be different.

Therefore, in this study, we assigned priority levels to trains which model different freights. We also considered safety and operational constraints that prevent the collision of trains as was done in \cite{rahman2014freight}. Proximity conflicts were considered for two trains travelling in the same
direction, and collision conflicts were considered for two trains travelling in opposite directions in the same segment.

With all the aforementioned constraints, a mathematical model with the notations described in Tables \ref{parameters_table}, \ref{Decision_variables}, and \ref{Binary_variables} was proposed.

\begin{table}[ht]
\caption{Binary variables of the scheduling model}
\label{Binary_variables}
\centering
\scriptsize
\resizebox{0.8\textwidth}{!}{%
\begin{tabular}{p{1.5cm}p{7.5cm}}
\hline
Symbol & Definition \\
\hline
$\alpha$ $_{pd}^{t\acute{t}}$ & =1 if train $t$ departs segment pd before train {t'}  \\
 & =0 otherwise\\[5pt]
$\beta$ $_{pdr}^{r\acute{r}}$ & =1 if train $r$ departs segment $pdr$ before train {r'}\\
 & =0 otherwise\\[5pt]
$\gamma$ $_{pd, pdr}^{t\acute{r}}$ & =1 if train $t$ departs segment pd before train $r$ \\ 
 & =0 otherwise
\end{tabular}%
}
\end{table}

Binary variables defined in Table \ref{Binary_variables} are considered to prevent safety conflicts, i.e., the first two variables are defined for the proximity conflicts and the last one for collision conflicts. The mathematical model of scheduling freight trains in a single line corridor is shown below:

\hspace*{-1cm}\vbox{\begin{alignat}{3}
& Min\sum\limits_{t \in T} {({\pi^t}(a_{SP(t)}^t - d_{EP(t)}^t} )) \label{eq1} \\
& s.t : \nonumber\\
& a_{pd}^t - d_{pd}^t \ge mr_{pd}^t &
\hspace{3cm},t \in T ; pd \in S \label{eq2}\\[5pt]
& a_{pdr}^r - d_{pdr}^r \ge mr_{pdr}^r   &,r \in T;pdr \in S \label{eq3}\\[5pt]
& d_{pd}^{t'} - a_{pd - 1}^t \ge Ul_{pd}^t +Lo_{pd}^t + dw_{pd}^t  &,t \in T,pd \in S \label{eq4}\\[5pt]
& d_{pdr}^{r'} - a_{pdr - 1}^r \ge Ul_{pdr}^r +Lo_{pdr}^r + dw_{pdr}^r &    ,r \in T,pdr \in S \label{eq5}\\[5pt]
& d_{pd}^{t'} - d_{pd}^t \ge ST_{pd} - M(1-\alpha_{pd}^{tt'})  &if \alpha_{pd}^{tt'} = 1 \label{eq6}\\[5pt]
& d_{pd}^{t} - d_{pd}^{t'} \ge ST_{pd} - M\alpha_{pd}^{tt'}   &if \alpha_{pd}^{tt'} = 0 \label{eq7}\\[5pt]
&d_{pdr}^{r'} - d_{pdr}^r \ge ST_{pdr} - M(1-\alpha_{pdr}^{rr'})  &if \beta_{pdr}^{rr'} = 1 \label{eq8}\\[5pt]
&d_{pdr}^{r} - d_{pdr}^{r'} \ge ST_{pdr} - M\beta_{pdr}^{rr'}   &if \beta_{pdr}^{rr'} = 0 \label{eq9}\\[5pt]
&\alpha _{pd}^{tt'} + \alpha _{pd}^{t't} \le1          &\forall t,t' \in T,pd \in S \label{eq10}\\[5pt]
&\beta _{pdr}^{rr'} + \beta _{pdr}^{r'r} \le1          &\forall r,r' \in T,pdr \in S \label{eq11}\\[5pt]
&\gamma _{pd,pdr}^{tr} + \gamma _{pd,pdr}^{rt}\le 1    &\forall t,r \in T,pd,pdr \in S,pd + pdr = n \label{eq12}\\[5pt]
&ar_{pdr}^r \le d_{pd}^t + M(1 - \gamma_{pd,pdr}^{rt}) &\forall t,r \in T,pd,pdr \in S,pd + pdr = n \label{eq13}\\[5pt]
&a_{pd}^t \le dr_{pdr}^r + M\gamma_{pd,pdr}^{rt}       &\forall t,r \in T,pd,pdr \in S,pd + pdr = n \label{eq14}
\end{alignat}}

As mentioned before, in this study, the objective function of scheduling freight trains was to minimize the total travel time of freight trains with respect to their priorities.

Trains must take a specific time to traverse the segments, which is defined by the train’s speed and the length of the segment (constraints \ref{eq2} and \ref{eq3} for departing and returning trains, respectively). Constraints \ref{eq4} and \ref{eq5} represent the time it takes for trains to dwell, load, and unload at all stations, and that their stop time cannot be less than the sum of these times.

Due to proximity conflicts, for two departing trains, the departure time for the following train must be more than the safe time after the departure time for the leading train (set of constraints \ref{eq6} and \ref{eq7} for departing trains, and constraints \ref{eq8} and \ref{eq9} for returning trains). On GAMS programming, each of the binary variables are shown two times in order to represent constraints of two following trains: constraints \ref{eq10} (for departing trains), \ref{eq11} (for returning trains), and \ref{eq12} (for all trains).

To ensure safe operation, the model is subject to more constraints. Safety constraints (\ref{eq13} and \ref{eq14}) illustrate that the departure time difference between two trains traversing in opposite directions must be greater than the travel time of the first departing train.

\subsection{The allocation problem} \label{allocation problem}

Goods were divided into two main categories: low-priority, and high-priority goods. We assigned priority to trains transporting high-priority goods in this section. We proposed a mathematical model to solve the problem of allocating the second category of goods, which may not consistently fill the trains’ capacity. We assumed that the first category of goods are in surplus; therefore, if the second category of goods cannot fill the trains’ capacity,the first category goods will also be loaded. However, if there are enough second category goods to fill the train's capacity, then only the second category of goods will be loaded. In addition, the following assumptions were made: \\

\begin{itemize}
    \item Freight weights may vary
    \item Trains’ weight capacity may vary
    \item Due to trains’ unique weight capacities and departure times, some freights may not be allocated to certain trains 
    \item At least 60 percent of each trains’ weight capacity is allocated to the second category of goods. (The whole weight capacity of trains can be allocated to the second category of freight)
    \item Freight has due dates at its destinations with penalties if tardy
    \item Freight has release dates at its start-stations.
\end{itemize}

In order to present the allocation model, we begin by stating the notation used in our model. Tables \ref{allocation_model_parameters} and \ref{Decision Variables_allocation} show the definition of sets and parameters, and decision variables used in the allocation model, respectively. Tardiness of freight $j$ is defined as the difference between the arrival time and the due date of freight $j$ as is shown below :
$$tardi_j=w{r_j} - {u_j} , t \in T;j \in J $$

\begin{table}[t]
\caption{Sets and parameters of the allocation model}
\label{allocation_model_parameters}
\centering
\scriptsize
\resizebox{0.7\textwidth}{!}{%
\begin{tabular}{p{1.5cm}p{6.5cm}}
\hline
Symbol & Definition \\
\hline
$j$ &  =\{1,2,3,4,5\}, set of freight \\ [5pt]
$w_j$ & Priority of freight $j$\ \\  [5pt]
$\delta$ $_j$ & Weight of freight $j$ \\ [5pt]
$\lambda$ $_t$ & Weight capacity of train $t$\\ [5pt]
${wr}_j$ &  Arrival time of freight $j$ to the end-station \\ [5pt]
$u_j$ & Due date of freight $j$ \\  [5pt]
${tardi}_j$ & Tardiness of freight $j$ \\ [5pt]
$O_j$ & Release date of freight $j$ \\  [5pt]
$M$ &  Sufficiently large constant 
\end{tabular}%
}
\end{table}

\begin{table}[t]
\caption{Decision Variables of the allocation model}
\label{Decision Variables_allocation}
\centering
\scriptsize
\resizebox{0.7\textwidth}{!}{%
\begin{tabular}{p{1.5cm}p{6.5cm}}
\hline
Symbol & Definition \\
\hline
$X_{jt}$ & =1 if freight $j$ allocates to train  $t$ \\
& =0 otherwise\\
\end{tabular}%
}
\end{table}

\begin{table}[b]
\caption{Solving time and efficiency of proposed formulation by GAMS software}
\label{solving_time}
\scriptsize
\resizebox{\textwidth}{!}{%
\begin{tabular}{lll}
\hline
Number of departing trains| & Number of returning trains| & Solving time of GAMS software \\
\hline
3 & 2 & 1 second \\
60 & 20 & 16 minutes \\
60 & 80 & 16 minutes \\
120 & 80 & 16 minutes
\end{tabular}}
\end{table}

The allocation model is shown below. The first sum of the objective function (shown in \ref{eq15}) minimizes the total penalty for tardy freight. In order to allocate higher priority freight to the scheduled trains, we subtracted the second sum from the objective function. Based on the freight’s weight and the weight capacity of the freight train, the second sum maximizes the allocation of higher priority freights to the scheduled trains.

\hspace*{-2cm}\vbox{\begin{alignat}{3}
&Min(\sum\limits_j {({w_j}.tardi_{j}) - \sum\limits_j {\sum\limits_t{({w_j}.{x_{j,t}}))} } } \label{eq15} \\
&s.t : \nonumber\\
&(0.6 \times {\lambda _t}) - M.{y_1} \le \sum\limits_j {({\delta _j}.{x_{j,t}}) \le } {\lambda _t}   &,t \in T;j \in J \label{eq16}\\[5pt]
&a_{SP}^t \le w{r_j} + M(1 - {x_{j,t}})   &,t \in T;j \in J \label{eq18}\\[5pt]
&\sum\limits_t {{x_{jt}} \le 1}   &,t \in T;j \in J \label{eq19}\\[5pt]
&d_{EP}^t \ge {O_j}.{x_{j,t}} + LO_{EP}^t  &,t \in T;j \in J \label{eq20}
\end{alignat}}

Constraint \ref{eq16} ensures that at least $60$ percent of the weight capacity of trains is allocated to the second category freight. If freight $j$ allocates to the train $t$ ($x_{j,t}=1$), then the arrival time of freight $j$ is equal to the arrival time of train $t$ at the end-station. Otherwise, $x_{j,t}=0$ , meaning that freight $j$ does not allocate to train $t$, and the arrival time of freight $j$ is equal to the arrival time of train $t$ to the end-station plus constant $M$ (constraint \ref{eq18}). Constraint \ref{eq19} ensures that each freight allocates to only one train. For all freight trains, if freight $j$ allocates to train $t$, then the departure time of the freight train must be greater than release date plus loading time of freight $j$ (constraint \ref{eq20}).

\section{Analysis of the mathematical models} \label{	Analysis of the mathematical models}
This section is organized to describe and illustrate the results of both aforementioned models in two different subsections. In this study, the mathematical model was solved for a single corridor with three segments (4 stations), three departing trains, and two returning trains. GAMS software solved the problem in 1 second. However, in order to evaluate the efficiency and solving time of the proposed formulation with GAMS software, we increased the number of trains in both directions. The results are shown in Table \ref{solving_time}.

To demonstrate the model, we generated data to illustrate the computational results. train’s priority was considered to be a random constant between zero and one. Additionally, loading, unloading, dwelling, and the safety time of trains were considered as random constants between zero and four hours. Due dates were considered equal to 10 hours for all trains and the release dates were developed as a random constant between time zero and five.

To illustrate the computational results of our proposed model, different figures are presented in this section. The bold green vertical lines at the end of each rectangle divide the segments. For all trains, the stop time is considered as the sum of loading, unloading, and dwelling time. At some stations, trains may not have loading, unloading or dwelling time. The black rectangles show the stop time of each train at each station. Two types of figures are illustrated in this section: train-time figures and train-location figures. In all figures, departing and returning trains are shown in blue and red, respectively. In all train-time figures, the vertical axis and horizontal axis are assigned to trains and time, respectively. Also, trains are assigned to the vertical axis and locations are assigned to the horizontal axis in all train-location figures. Table {\ref{colors}} illustrates the definition of colors used in the results figures.

\begin{table}[t]
\caption{Colors used in the Figures}
\label{colors}
\centering
\scriptsize
\resizebox{0.7\textwidth}{!}{%
\begin{tabular}{p{1.5cm}p{7.5cm}}
\hline
Color & Definition of colors \\
\hline
\parbox[c]{1em}{\includegraphics[height=8mm]{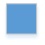}} & Traversing time of departing trains at segments \\
\parbox[c]{1em}{\includegraphics[height=8mm]{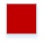}} & Traversing time of returning trains at segments \\
\parbox[c]{1em}{\includegraphics[height=8mm]{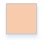}} & Loading time of departing trains at stations \\
\parbox[c]{1em}{\includegraphics[height=8mm]{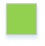}} & Unloading time of departing trains at stations \\
\parbox[c]{1em}{\includegraphics[height=8mm]{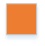}} & Dwelling time of departing trains at stations \\
\parbox[c]{1em}{\includegraphics[height=8mm]{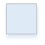}} & Loading time of returning trains at stations \\
\parbox[c]{1em}{\includegraphics[height=8mm]{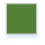}} & Unloading time of returning trains at stations \\
\parbox[c]{1em}{\includegraphics[height=8mm]{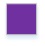}} & Dwelling time of returning trains at stations 
\end{tabular}%
}
\end{table}

\subsection{Computational results of the scheduling problem}
This section illustrates the computational results of the scheduling problem.
\subsubsection{Scheduling single line with one direction}\

In order to evaluate the proposed model and its proximity constraints, we defined three different objective functions for single lines with one direction: minimizing the total departure time for all trains from the start-station, minimiz\-ing the total arrival times for all trains to the end-station, and minimizing the total travel time of all trains. 

Based on the aforementioned notations, the first objective function for evaluating the proposed model in a single line with one direction was formulated as $Min\sum_{t} d_{EP}^t$ and Figure \ref{fig_1a} shows the results. 

As illustrated in Figure \ref{fig_1a}, train one departs from the first station at time 1 and arrives to the end of segment 1, i.e., $pd1$ , at time 3. It spends 3 hours loading, unloading, and dwelling at station 2 (beginning of segment 2) and arrives to the last station at time 13. Because the goal is to minimize the sum of the trains’ departure times, the model reduces the speed of trains two and three at segment one. The objective was minimized to three, which is the least objective value that can be achieved, to the best of our knowledge.

\begin{figure*}[t!]
    \centering
    \begin{subfigure}[t]{0.5\textwidth}
        \centering
        \includegraphics[height=1.5in]{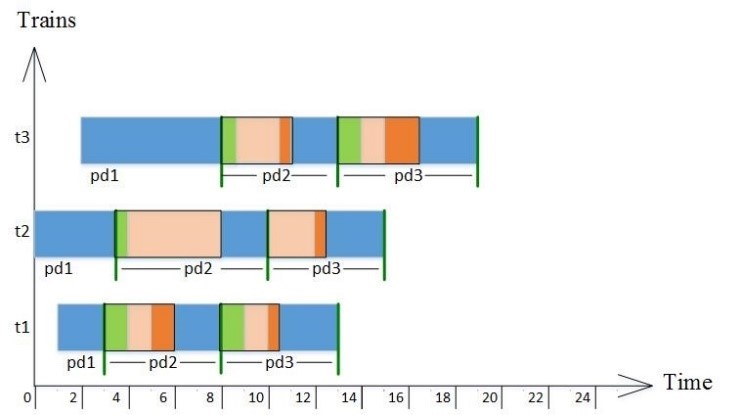}
        \caption{Minimizing total departure time from start-station}
        \label{fig_1a}
    \end{subfigure}%
    ~ 
    \begin{subfigure}[t]{0.5\textwidth}
        \centering
        \includegraphics[height=1.5in]{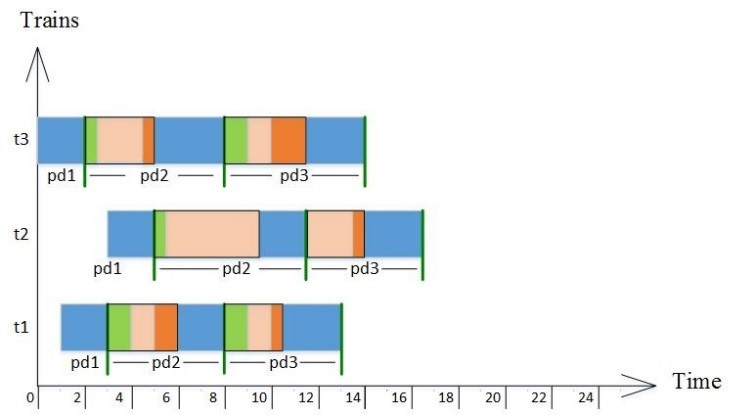}
        \caption{Minimizing total arrival time to the end-station}
        \label{fig_1b}
    \end{subfigure}
     ~ 
    \begin{subfigure}[t]{0.5\textwidth}
        \centering
        \includegraphics[height=1.5in]{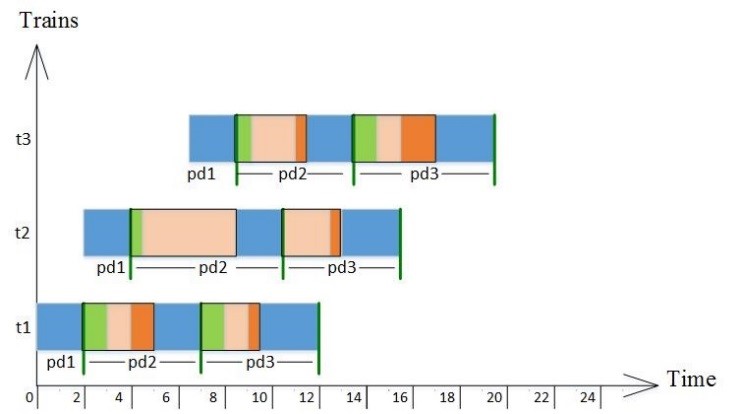}
        \caption{Minimizing total travel time}
        \label{fig_1c}
    \end{subfigure}
    \caption{Gantt charts of scheduling freight trains on a single line with one direction}
\end{figure*}

The second objective function is shown as $Min\sum_{t} a_{SP}^t$ and the computational results for three trains and three segments are shown in Figure \ref{fig_1b}. As shown in Figure \ref{fig_1b}, the total arrival time for all trains is 43.5, the minimum objective value.

The third objective function is minimizing the total travel time of trains, shown in Figure \ref{fig_1c}. As expected from the objective function, the model minimizes the travel times of trains; the final value is 38.5, which is the least value for the objective. Thus, based on all three evaluation models, the constraints of proximity conflicts gave us promising objective values.
 
 \begin{figure}[t]
     \centering
     \includegraphics[width=0.7\textwidth]{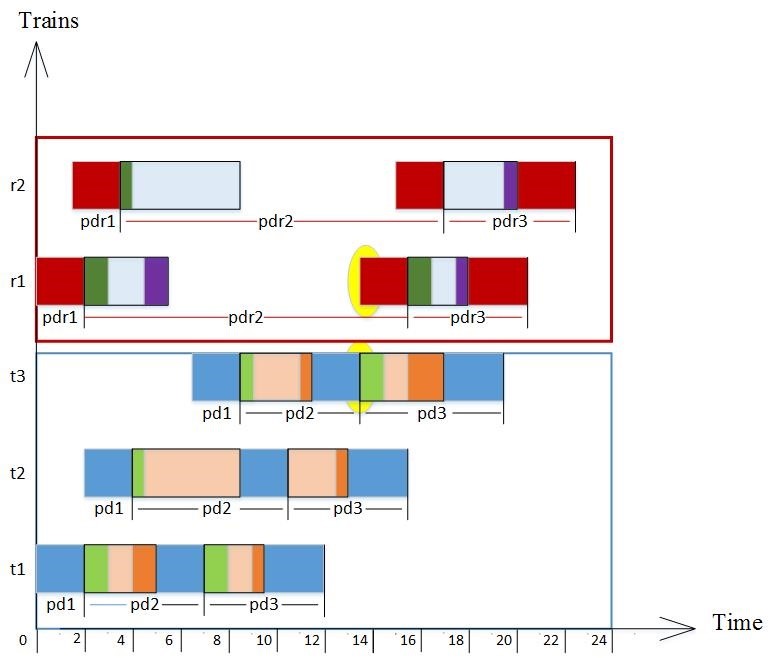}
     \caption{Gantt charts of scheduling model on a single line corridor (Train-Time)}
     \label{fig_2}
 \end{figure}
 
\subsubsection{Scheduling single line corridor (two directions)}\

we illustrate the results of the complete scheduling model. We considered two names for the segments, one for departing trains and one for returning trains. Figure \ref{fig_2} illustrates the Gantt chart of the scheduling model on a single line with three segments. departing trains traverse segments 1, 2 and 3, which are shown as $pd1$, $pd2$ and $pd3$, respectively. Returning trains first traverse segment 3 and then segments 2 and 1, shown as $pdr1$, $pdr2$ and $pdr3$. In other words, segment 1 is the same for both directions of trains; the only difference is the naming of them for departing and returning trains, i.e., segment 1 is $pd1$ for departing trains and $pdr3$ for returning trains. Computational results for an example of three departing and two returning trains are shown in Figure \ref{fig_2}. The upper red-bordered box outlines the returning trains, while the lower blue-bordered box shows the departing trains. Departing trains traverse segments without any interruption, while the returning trains are interrupted for the safety of departing trains, which is shown by the yellow ellipse. This was due to the safety constrains combined with the fact that departing trains were given higher priority.

The train-location graph is shown below. Table \ref{table_arrows} describes the symbols used in the train-location figure.

\begin{table}[t]
\caption{Arrows and their definitions used in the train-location figure}
\label{table_arrows}
\centering
\scriptsize
\resizebox{0.7\textwidth}{!}{%
\begin{tabular}{p{1.5cm}p{6.5cm}}
\hline
Symbols & Definition \\
\hline
\parbox[c]{1em}{\includegraphics[height=4mm]{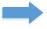}} & Departure of departing trains from a segment \\[5pt]
\parbox[c]{1em}{\includegraphics[height=4mm]{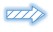}} & Arrival of departing trains to a segment \\[5pt]
\parbox[c]{1em}{\includegraphics[height=4mm]{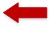}} & Departure of returning trains from a segment \\[5pt]
\parbox[c]{1em}{\includegraphics[height=4mm]{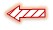}} & Arrival of returning trains to a segment
\end{tabular}%
}
\end{table}

 \begin{figure}[ht]
     \centering
     \includegraphics[width=0.8\textwidth]{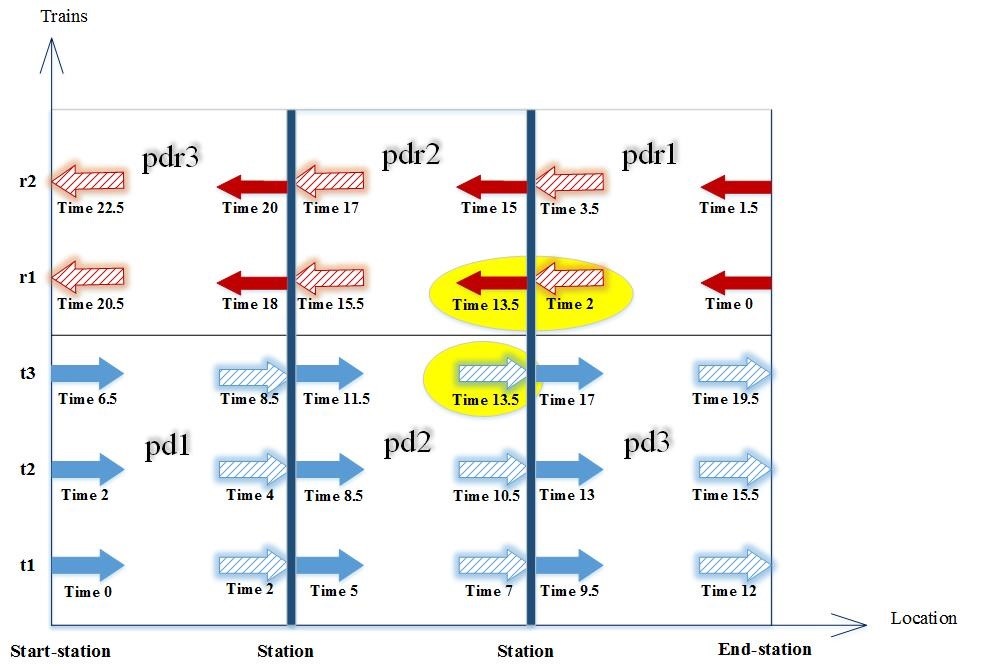}
     \caption{Departing and returning trains scheduling on a single line corridor}
     \label{fig_3}
 \end{figure}

The vertical rectangles shown in Figure \ref{fig_3} illustrate the segments, i.e., segment 1 for departing trains is shown as $pd1$ and as $pdr3$ for returning trains. The departure and arrival times of each train to each segment are shown under the arrows. Train $t3$  departures from segment 2 ($pd2$) at time 11.5 and arrives to the next station at time 13.5. Due to the safety constraints, two facing trains cannot simultaneously traverse the same segment. Thus, trains $r1$  and $r2$ must stop until train $t3$ traverses the second segment. Therefore, after time 13.5 (when train $t3$ finishes traversing segment 2), train $r1$ can start its traversing at segment 2. In order to highlight these interruptions, similar to Figure \ref{fig_2}, the interruptions are shown using yellow ellipses.

\begin{table}[ht]
\caption{Allocation of freights to trains}
\label{freights_allocation}
\centering
\scriptsize
\resizebox{0.8\textwidth}{!}{%
\begin{tabular}{p{1.5cm}p{1.5cm}p{1.5cm}p{1.5cm}p{1.5cm}p{1.5cm}l}
\hline
$x_{j,t}$ & $j_1$ & $j_2$ & $j_3$ & $j_4$ & $j_5$ \\[5pt]
\hline
$t_1$ & 0 & 0 & 0 & 1 & 0 \\[5pt]
$t_2$ & 1 & 1 & 1 & 0 & 0 \\[5pt]
$t_3$ & 0 & 0 & 0 & 0 & 1 
\end{tabular}%
}
\end{table}

\begin{table}[ht]
\caption{Outputs of the allocation model}
\label{allocation_output}
\centering
\scriptsize
\resizebox{0.8\textwidth}{!}{%
\begin{tabular}{p{1.5cm}p{1.5cm}p{1.5cm}p{1.5cm}p{1.5cm}p{1.5cm}l}
\hline
Freight & $j_1$ & $j_2$ & $j_3$ & $j_4$ & $j_5$ \\[5pt]
\hline
$wr_j$ & 18.5 & 18.5 & 18.5 & 12 & 23.5 \\ [5pt]
$tardi_j$ & 8.5 & 8.5 & 8.5 & 2 & 13.5 \\  [5pt]
\end{tabular}%
}
\end{table}

\subsection{Computational results of the allocation model} The computational results of the allocation model have three output variables: the allocation of freight to trains ($x_{j,t}$), the arrival time of freight to the end-station (${wr}_{j}$), and the tardiness of freight at the end-station. Five different freight loads were solved with the model and the outputs are shown in (Table \ref{allocation_output}).

Table \ref{freights_allocation} displays that, based on the freight weights, weight capacity of trains, and due dates of freights, freights 1, 2, and 3 were allocated to train 2, and freights 4 and 5 were allocated to trains 1 and 3, respectively. As mentioned in the previous section, we wanted to allocate each freight load to a single train, and as the results show in Table \ref{freights_allocation}, this constraint is also satisfied.

\section{Heuristic Algorithm}
As real-world railway systems have constraints that do not easily fit into a simple mathematical formulation in real large-scale problems \cite{cai1998greedy}, a heuristic algorithm was used to address the real train scheduling. It makes it possible to solve many of such problems and provides alternatives that seem to be the best at that moment. The proposed heuristic can not only deal with all constraints of the mathematical model, but can also match with other situations flexibly. The algorithm framework generates feasible solutions according to all constraints based on the freight priority ($w_j$) and average speed of each train and then it inspires the PSO to enhance the solution quality; this way, each solution is referred to two personal and global best solutions. The algorithm selects the trains that have had the worst effect on the fitness function and tries to find a new improved sequence for the trains.

As mentioned earlier, both the allocation and scheduling models have inherent complexities and can be solved only separately for small problems. Since scheduling is not applicable without an allocation program, it is not practically possible to expect the above models to be individually trouble-shooters for the industry. It should be admitted, of course, that this can be partially overcome by such methods as the benders decomposition algorithm. In this research, however, we tried to develop a heuristic algorithm to obtain an acceptable solution within the constraints of both the scheduling and allocation problems.

The algorithm to be presented next has been somehow inspired by the particle swarm algorithm. In the latter, there is a limitation that each solution is consider\-ed as a particle in an n-dimensional space, and the improving movement of each particle toward the reference particles (personal, local, and global) occur integrated across all the particle dimensions. But in the proposed algorithm, each dimension of the particle can be changed and improved independently (of other dimensions) which helps the solution to remain unchanged in the appropriate dimensions. First, use is made of a reproduction algorithm to generate the required number of initial solutions (Algorithm I) and then a reference set is selected to improve each solution based on the guidelines explained in Algorithm \rom{2}.

\begin{table}[]
\caption*{Algorithm \rom{1}, produce the primary generation}
\label{Algorithm I }
\centering
\footnotesize
\begin{tabular}{p{1cm}p{12cm}}
\hline
\multicolumn{2}{l}{do for $n$ solutions: (outputs in each iteration: $\zeta_e, {\zeta_e^\prime}, F_{\zeta_e}$ and $F_{\zeta_e}^\prime)$} \\
i.   & $k = k^{\prime} = 1$ \\
     & $\zeta_e = \zeta_e^\prime = \{\}$ \\
ii.  & Arrange the trains based on $(w_j*\delta_j)/\lambda_j$  in an ascending order. The set of departing trains $\rightarrow \Psi_t$.The set of returning trains $\rightarrow \Psi_t^\prime$. \\
iii. & Calculate the selection probability of each train. \\
iv.  & Using Monte-Carlo method select trains $st_k$ and $st_{k^\prime}$ from sets $\Psi_t$ and $\Psi_t^\prime$, and update ordered sets $\zeta_e$ and $\zeta_e^\prime$. \\
     & $\zeta_e=\zeta_e +{st_k}$\\
     & ${\zeta_e^\prime}={\zeta_e^\prime} +st_{k^\prime}$\\
     & $\psi_e=\psi_e -st_{k^\prime}$\\
     & ${\psi_e^\prime} = {\psi_e^\prime} - st_{k^\prime}$\\
v.   & For every train in sets $\zeta_e$ and $\zeta_e^\prime$, calculate the earliest arrival ($f_k$) considering all constraints in both lines and create sets  $F_{\zeta_e}$ and $F_{\zeta_e^\prime}$. \\ 
\hline
\end{tabular}
\end{table}

Another difference between this algorithm and the particle swarm algorithm is in selecting the reference particle. In the latter, three references (one personal, one local and one global) are selected for each particle in any iteration and the improving movement vector is obtained from the resultant of the movement of each particle toward the three mentioned references. In the proposed algorithm, the improving movement is inspired by only one solution, but since the mentioned reference is random, it results in a escape from the local optimum trap. For each step, $\rho$ percent of its best solutions is considered and one of them is randomly select as the reference particle of that step. It is worth mentioning that, based on the $\nu_i$ index, only some trains (dimensions) are nominated for improvement, which has a significant impact on achieving good quality solutions.

\begin{table}[]
\caption*{Algorithm \rom{2}, heuristic improvement}
\label{Algorithm \rom{2}}
\centering
\footnotesize
\begin{tabular}{p{1cm}p{12cm}}
\hline
\multicolumn{2}{l}{Using Algorithm \rom{1}, produce set $\omega$ including $n$ feasible solutions and $n$ sets of $\zeta_e, \zeta_e^\prime, F_{\zeta_e}$ and $F_{\zeta_e^\prime}$} \\
\multicolumn{2}{l}{$l=1$} \\
\multicolumn{2}{l}{Do the following steps based on the termination condition.} \\
i.   & For each member of  $\Omega$, calculate the fitness function and then determine the best ($\phi_l^b$) and the worst ($\phi_l^w$)\\
ii.  & Opt $\rho$ percent of the high quality solutions of set $\Omega \rightarrow \omega_\rho$ \\
iii. & Select a member of  $\omega_\rho$ randomly as the reference of $\Omega_i \rightarrow {\omega_\rho}^i$ \\
iv.  & Calculate index $\nu_i$  as the number of corrections in the $\Omega_i$ sequences.\\
     & $\nu_i = \alpha((\phi_l^b - \phi_l^w)/(\phi_{\Omega_i} - \phi_l^w))$, $0<\alpha<1$\\
v.   & Based on $f_k$ and $f_(k+1)$, select $\nu_i$ trains that have the highest effect on the bad quality of the fitness function.\\ 
vi.  & Update $\zeta_e$ and $\zeta_e^\prime$; the sequences of $\nu_i$ selected trains will be reordered according to its reference ($\omega_\rho^i$)\\
vii. & Use the OPT algorithm to achieve a better solution based on the updated $\zeta_e$ and $\zeta_e^\prime$, update $F_{\zeta_e}$ and $F_{\zeta_e^\prime}$\\
\hline
\end{tabular}
\end{table}

\begin{table}[]
\caption{Updating a 10-train solution based on the references}
\label{10train}
\centering
\scriptsize
\begin{tabular}{p{2cm}p{0.7cm}p{0.7cm}p{0.7cm}p{0.7cm}p{0.7cm}p{0.7cm}p{0.7cm}p{0.7cm}p{0.7cm}p{0.7cm}}
\hline
 &  & * &  & * &  &  &  &  &  & * \\
$\Omega_i$   & 1 & 2 & 3 & 4 & 5 & 6 & 7 & 8 & 9 & 10 \\
$\omega_p^i$ & 4 & 3 & 5 & 6 & 7 & 8 & 10 & 1 & 2 & 9 \\
\hline
Updated $\Omega_i$ & 10 & 1 & 4 & 3 & 5 & 6 & 7 & 8 & 2 & 9 \\ 
\end{tabular}
\end{table}

To understand step 6 better, Table \ref{10train} shows how each solution is updated based on the reference, OPT is a well-known local search algorithm which improves a solution by changing the sequence of trains. The pseudo code related to the algorithm is as follows:

\begin{table}[]
\caption*{The OPT algorithm}
\label{OPT algorithm}
\centering
\begin{tabular}{p{1cm}p{8cm}}
\hline
\multicolumn{2}{l}{For each member of $\Omega_i$ except the last one:} \\
& Change the sequence of the considered train to the next one \\
& Calculate the new $\omega_{\Omega_i}$ \\
& If new $\omega_{\Omega_i} < old \omega_{\Omega_i}$ then update $\Omega_i$\\
\hline
\end{tabular}
\end{table}

\begin{table}[ht]
\caption{Comparison of the exact and heuristic algorithms}
\label{Comparison of teo algorithms}
\centering
\scriptsize
\begin{tabular}{ccccccccc}
\hline
\multirow{2}{*}{trains} & \multicolumn{3}{c}{Heuristic algorithm} & \multicolumn{1}{l}{} & \multicolumn{4}{c}{Exact algorithm} \\ \cline{2-4} \cline{6-9} 
 & \multicolumn{1}{l}{CPU time} & \multicolumn{1}{l}{Gap\%} & \multicolumn{1}{l}{Fitness function} & \multicolumn{1}{l}{} & \multicolumn{1}{l}{CPU time} & \multicolumn{1}{l}{Gap\%} & \multicolumn{1}{l}{Lower bound} & \multicolumn{1}{l}{Upper bound} \\ \hline
6 & 1.7 & 0 & 152 &  & 1 & 0 & 152 & 152 \\
20 & 64 & 0.078 & 781 &  & 52 & 0.05 & 724 & 762 \\
40 & 512 & 0.141 & 1617 &  & 1500 & 0.191 & 1416 & 1687 \\
60 & 1728 & 0.195 & 2414 &  & 5000 & 0.38 & 2019 & 2789 \\
80 & 4096 & 0.249 & 3621 &  & 5000 & 0.45 & 2897 & 4015 \\
100 & 8000 & - & 4856 &  & - & - & - & - \\ 
\end{tabular}
\end{table}

To test the algorithm, some sets of 6, 20, 40, 60, 80 and 100 trains traveling to and fro in one line with different speeds were generated and 10 different freights with different weights and priorities for loading and sending were considered. These sets were solved using exact and heuristic methods (the former covered only the scheduling not the allocation) and the results were compared (Table \ref{Comparison of teo algorithms}). To analyze the solutions, it is important to note that the gap found in the heuristic solution has been obtained based on the lower bound of the exact solution, and the reason why this gap does not seem very appropriate is that this lower bound is considerably lower than the optimal value. As an example, for the 80-train set, the gap related to the heuristic algorithm has been obtained using 2897 (Table \ref{Comparison of teo algorithms}), which is less than the optimal value and makes it inappropriate.

An important point in the proposed algorithm is the use of OPT algorithm, which greatly affects the solution improvement. Figure \ref{fig_4} shows two solutions for the 60-train set with and without using OPT. Although the latter does not have a significant effect on the solution improvement in early stages, and the main burden of the search engine is on the original algorithm, from step 280 onward, it shows its efficiency and creates a meaningful gap between the two solutions. As shown, both algorithms have the required convergence in final iterations, but the improved one converges in a better orbit.

 \begin{figure}[H]
     \centering
     \includegraphics[width=0.8\textwidth]{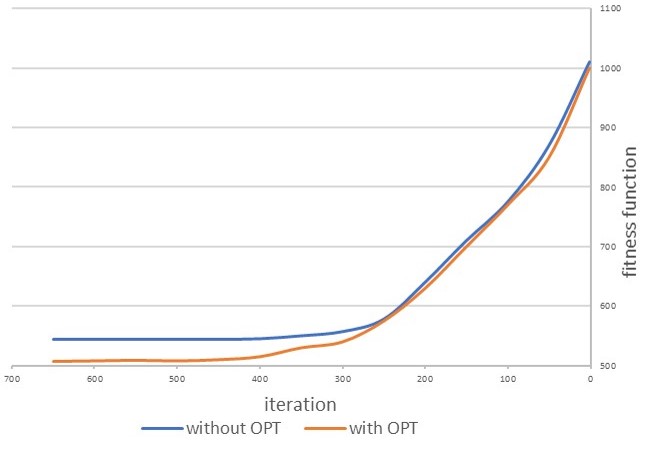}
     \caption{Efficiency of the OPT algorithm}
     \label{fig_4}
 \end{figure}

\section{Conclusion}
This paper addressed the scheduling and freight allocation to freight trains on single-line, to and fro routes. First, both issues were modeled separately and analyzed thoroughly, and then a heuristic algorithm was proposed. To this end, use was made of a PSO-inspired heuristic because the models were complex. Results showed that the proposed algorithm performed quite well and the use of a daemon algorithm, called OPT, gave it even a better performance. The notable point in this study, besides proposing an appropriate algorithm for scheduling and freight allocation to trains, is using the concept of the PSO algorithm and presenting a new concept of “change” that can be applied to other problems (e.g. routing). In this concept, each particle takes a shape more similar to the reference particle at each step rather than moving toward it as in the PSO algorithm.
Accordingly, it can be claimed that this paper proposes two different viewpoints for future studies: First, developing the mathematical model and then the solution algorithm considering the facts and scenarios of the scheduling problems and second, using the framework of the proposed algorithm to solve such other problems as the “travel salesman” or the VRP vehicle routing.

}
\end{flushleft}

\bibliographystyle{IEEEtran}
\bibliography{RefUserAuth}

\end{document}